# Computing special values of motivic $L$-functions

Tim Dokchitser

University of Durham, United Kingdom

**Abstract**

We present an algorithm to compute values $L(s)$ and derivatives $L^{(k)}(s)$ of $L$-functions of motivic origin numerically to required accuracy. Specifically, the method applies to any $L$-series whose $\Gamma$-factor is of the form $A^s \prod_{i=1}^{d} \Gamma(\frac{s+\lambda_j}{2})$ with $d$ arbitrary and complex $\lambda_j$, not necessarily distinct. The algorithm relies on the known (or conjectural) functional equation for $L(s)$.

## 1 Introduction

Many $L$-series in number theory and algebraic geometry can be interpreted as $L$-series of motives over number fields. For instance, Riemann and Dedekind $\zeta$-function, Dirichlet and Artin $L$-series and $L$-series of elliptic curves are of this kind. These are all of the form $L(X, V, s)$ associated to $V = H^i(X)$ or a "motivic" subspace $V \subset H^i(X)$ of a projective algebraic variety $X/K$.

Given such series

$$L(s) = \sum_{n=1}^{\infty} \frac{a_n}{n^s}, \qquad \operatorname{Re} s \gg 0, \qquad (1)$$

standard conjectures state that $L(s)$ extends to a meromorphic function on the whole of $\mathbf{C}$ and satisfies a functional equation of a predicted form. The Riemann hypothesis tells where the zeroes of $L(s)$ are located and, finally, various conjectures relate values of $L(s)$ at integers to arithmetic invariants of $X$. For instance, the Birch-Swinnerton-Dyer [2], Zagier [26], Deligne-Beilinson-Scholl [1, 20] and Bloch-Kato [3] conjectures are examples of these.

While the aforementioned conjectures remain unproved in the vast majority of cases, a lot of work has been done to provide numerical evidence for some of them in low-dimensional cases. This applies especially to the Riemann hypothesis for the Riemann $\zeta$-function [24], Dirichlet and Artin $L$-series [6, 13, 15, 19, 23] and $L$-series $L(E, H^1, s)$ of elliptic curves [9]. Another well-studied case is the Birch-Swinnerton-Dyer conjecture [2, 5] for $L(E, H^1, s)|_{s=1}$ where $E/\mathbf{Q}$ is an elliptic curve.

To perform this kind of calculations one needs an efficient algorithm to compute $L(s)$ (or, more precisely, its analytic continuation) numerically to required precision for a given complex $s$. Such algorithms are usually based on writing $L(s)$ as a series in special functions associated to the inverse Mellin transform of the $\Gamma$-factor of $L(s)$. In the cases mentioned above these special functions are incomplete Gamma functions for $\dim V = 1$ (Riemann $\zeta$-function, Dirichlet characters) and incomplete Bessel functions for $\dim V = 2$ (modular forms, elliptic curves).

The original motivation for this paper was to study evidence for the "standard conjectures" in case of higher-dimensional motives, for instance curves of genus $g > 1$ etc. One then needs an efficient method to compute $L(X, V, s)$ in cases $\dim V > 2$ as well. In this paper we present such a method. Namely, for an arbitrary motivic $L$-series for which meromorphic continuation and the functional equation are assumed, the algorithm numerically verifies the functional equation and allows to compute the values $L(s)$ and derivatives $L^{(k)}(s)$ for arbitrary complex $s$ to required precision.

The scheme presented here has been implemented as a PARI script and is available on [7]. This includes examples of computations with the Riemann $\zeta(s)$, Dirichlet $L$-functions, Dedekind $\zeta$-function, Shintani's $\zeta$-function, $L$-series of modular forms and those associated to curves $C/\mathbf{Q}$ of genus 1,2,3 and 4. Also note that the formulae described in the paper can be used in any other environment as long as it provides arbitrary precision arithmetic, complex numbers, Laurent series and the Taylor series expansion of the $\Gamma$-function.

The structure of the paper is as follows: in §2 we start with generalities on the invariants of $L$-functions and outline the algorithm. The approach used here is standard and has been used in most of the algorithms to compute $L$-functions (e.g. [15, 19, 23, 24]). In §3 and §4 we deduce power series expansions of general Meijer G-functions required in the computations. In §5 we present asymptotic expansions at infinity of the same functions and continued fraction expansions associated to those. Then §6 summarises the algorithm and addresses implementation and accuracy issues. Finally, §7 contains some remarks on working with $L$-functions for which not all of the invariants are known.

The author is extremely grateful to Don Zagier for suggesting to work on this project to begin with and also for numerous explanations, ideas and suggestions which have finally lead to this work. Just saying that the algorithm and this paper would not exist without his influence illustrates only a small fraction of the truth. The author would also like to thank the stimulating atmosphere of the Max-Planck Institut für Mathematik in Bonn where most of this work has been carried out.

## 2 Motivic $L$-functions

Suppose we are given an $L$-series,

$$L(s) = \sum_{n=1}^{\infty} \frac{a_n}{n^s}, \qquad a_n \in \mathbf{C}.$$

We make the following three assumptions on $L(s)$:

**Assumption 2.1.** The coefficients of $L(s)$ grow at most polynomially in $n$, that is $a_n = \mathcal{O}(n^\alpha)$ for some $\alpha > 0$. Equivalently the defining series for $L(s)$ converges for $\text{Re } s \gg 0$.

**Assumption 2.2.** $L(s)$ admits a meromorphic continuation to the entire complex plane. There exist: *weight* $w \geq 0$, *sign* $\epsilon \in \mathbf{C}$, real positive *exponential factor* $A$ and the *$\Gamma$-factor*

$$\gamma(s) = \Gamma\Big(\frac{s+\lambda_1}{2}\Big) \cdots \Gamma\Big(\frac{s+\lambda_d}{2}\Big)$$

of *dimension* $d \geq 1$ and *Hodge numbers* $\lambda_1, \ldots \lambda_d \in \mathbf{C}$, such that

$$L^*(s) = A^s \, \gamma(s) \, L(s)$$

satisfies a *functional equation*[1]

$$L^*(s) = \epsilon \, L^*(w-s) \,. \tag{2}$$

**Assumption 2.3.** $L^*(s)$ has finitely many simple poles $p_j$ with residues $r_j = \mathrm{res}_{s=p_j} L^*(s)$ and no other singularities.

**Remark.** Even general motivic L-functions have much more specific parameters. Usually $a_n$ lie in the ring of integers of a fixed number field (most often $\mathbf{Z}$), $A = \sqrt{N}/\pi^{d/2}$ (with *conductor* $N \in \mathbf{Z}$), $\lambda_k$ are integers (or even $\lambda_k \in \{0,1\}$) and $\epsilon = \pm 1$. Moreover, $L^*(s)$ is usually entire and there is a product formula for $L(s)$. However, these additional assumptions do not simplify the algorithm and there are some L-functions not of motivic origin (e.g. Shintani's $\zeta$-function [22]), to which the algorithm applies. So we do not require more than stated above. The last assumption that the poles of $L^*(s)$ are simple is not essential either (see the discussion below).

**Example 2.4.** The following table contains some well-known examples of L-series satisfying our assumptions and their basic invariants. For every one of those L-functions, the exponential factor is of the form $A = \sqrt{N}/\pi^{d/2}$ with $N \in \mathbf{Z}$.

| $L(s)$ | Description | $w$ | $d$ | $(\lambda_j)$ | $N$ | $\epsilon$ | $(p_j)$ |
|---|---|---|---|---|---|---|---|
| $\zeta(s)$ | Riemann $\zeta$-function | 1 | 1 | (0) | 1 | 1 | (0,1) |
| $L(\chi,s)$ $L(\bar\chi,s)$ | $\chi$ primitive Dirichlet character mod $N$ | 1 | 1 | (0), $\chi(-1)=1$ (1), $\chi(-1)=-1$ | $N$ | $\|\epsilon\|=1$ | |
| $\zeta(F,s)$ | Dedekind $\zeta$-function $[F:\mathbf{Q}]=d$ | 1 | $d$ | $(0,\ldots,0,1,\ldots,1)$ $d-\sigma, \sigma$ times | $\|\Delta_F\|$ | 1 | (0,1) |
| $L(f,s)$ | $f$ modular form of weight $k$ on $\mathrm{SL}_2(\mathbf{Z})$ | k | 2 | (0, 1) | 1 | $(-1)^k$ | (0,k) |
| $L(f,s)$ | $f$ cusp form of weight $k$ on $\mathrm{SL}_2(\mathbf{Z})$ | k | 2 | (0, 1) | 1 | $(-1)^k$ | |
| $L(f,s)$ | $f$ Hecke cusp form of weight $k$ on $\Gamma_0(N)$ | k | 2 | (0, 1) | $N$ | $\pm 1$ | |
| $L(E,s)$ | $E/\mathbf{Q}$ elliptic curve of conductor $N$ | 2 | 2 | (0, 1) | $N$ | $\pm 1$ | |
| $L(C,s)$ | $C/\mathbf{Q}$ genus $g$ curve of conductor $N$ | 2 | $2g$ | $(0,\ldots,0,1,\ldots,1)$ $g, g$ times | $N$ | $\pm 1$ | |
| $\zeta_{\mathrm{Sh}}(s)$ | Shintani's $\zeta$-function | 1 | 4 | $(0, 1, \frac{1}{6}, -\frac{1}{6})$ | $2^4 3^3$ | 1 | $(0, \frac{1}{6}, \frac{5}{6}, 1)$ |

---
[1] Functional equation may also involve two different L-functions, see Remark 2.5

In the second row $L(\chi, s)$ satisfies a functional equation which involves the "dual" $L$-function associated to the complex conjugate character $L(\bar{\chi}, s)$ (see Remark 2.5 below). In the third row $\Delta_F$ is the discriminant of $F/\mathbf{Q}$ and $\sigma$ is the number of pairs of complex embeddings.

For the last (non-motivic) example, see Shintani's original paper [22]. For the rest (and more), see [18], Chapter 4 and articles in [12] for references and additional information. For actual $L$-series computations in the above cases, see [7].

Given an $L$-function which satisfies 2.1–2.3 we would like to

(a) Give a numerical verification of the functional equation for $L(s)$.

(b) For a given $s_0 \in \mathbf{C}$ and an integer $k \geq 0$ determine the $k$-th derivative $L^{(k)}(s_0)$ to given precision.

To this end define $\phi(t)$ to be the inverse Mellin transform of $\gamma(s)$, that is

$$\gamma(s) = \int_0^\infty \phi(t)\, t^s \, \frac{dt}{t} \,. \tag{3}$$

Henceforth we let $s$ denote a complex number and $t$ a positive real. The function $\phi(t)$ exists (for real $t > 0$ that is) and it decays exponentially for large $t$ (see §3). In particular, the following sum converges exponentially fast,

$$\Theta(t) = \sum_{n=1}^\infty a_n \, \phi\!\left(\frac{nt}{A}\right) \tag{4}$$

This function is defined so that $L^*(s)$ becomes the Mellin transform of $\Theta(t)$,

$$\int_0^\infty \Theta(t) t^s \frac{dt}{t} = \int_0^\infty \sum_{n=1}^\infty a_n\, \phi\!\left(\frac{nt}{A}\right) t^s \frac{dt}{t} = \sum_{n=1}^\infty a_n \int_0^\infty \phi\!\left(\frac{nt}{A}\right) t^s \frac{dt}{t}$$
$$= \sum_{n=1}^\infty a_n \int_0^\infty \phi(t) \left(\frac{At}{n}\right)^s \frac{dt}{t} = A^s \sum_{n=1}^\infty \frac{a_n}{n^s}\, \gamma(s) = L^*(s)\,. \tag{5}$$

By Mellin's inversion formula

$$\Theta(t) = \int_{c-i\infty}^{c+i\infty} L^*(s) t^{-s} ds, \qquad \operatorname{Re} c \gg 0,$$

if $c \in \mathbf{C}$ is chosen to lie to the right of the poles of $L^*(s)$. By the assumed functional equation (2) for $L^*(s)$,

$$\Theta(1/t) = \int_{c-i\infty}^{c+i\infty} L^*(s) t^s ds = t^w \int_{c-i\infty}^{c+i\infty} \epsilon L^*(w-s) t^{s-w} ds$$
$$= t^w \epsilon \int_{w-c-i\infty}^{w-c+i\infty} L^*(s) t^{-s} ds$$

which is almost an expression for $\epsilon\, t^w \Theta(t)$, except that the integration path lies on the *left* of the poles of $L^*(s)$. Shifting this path to the right we pick up residues of $L^*(s)t^{-s}$ at every pole of $L^*(s)$. Consequently $\Theta(t)$ enjoys a functional equation,

$$\Theta(1/t) \;=\; \epsilon\, t^w \Theta(t) - \sum_j r_j t^{p_j}\,. \qquad (6)$$

Note that the assumption that $L^*(s)$ has simple poles in not essential. If the poles are of higher order, the residues of $L^*(s)t^{-s}$ also involve some $\log t$-terms, so (6) and (10) below have extra terms but this does not affect the reasoning elsewhere.

In §3 and §5 we describe how to compute $\phi(t)$ for $t>0$ for a given $\Gamma$-factor $\gamma(s)$. Then $\Theta(t)$ can be also effectively computed numerically since (4) converges exponentially fast.

Now we can answer the first question, that of numerical verification of the functional equation for $L^*(s)$. Pick $t>0$ and check that (6) holds numerically for this $t$. In fact (6) holds for *all* $t$ if and only if the functional equation (2) is satisfied. Note that having such a verification is useful when not all of the invariants of $L(s)$ are known (see §6).

**Example.** Let $L(s) = \zeta(s) = \sum_{n=1}^{\infty} n^{-s}$ be the Riemann $\zeta$-function. Then

$$a_n \equiv 1, \quad w=1, \quad \epsilon=1, \quad A=\frac{1}{\sqrt{\pi}}, \quad d=1, \quad \gamma(s)=\Gamma(\tfrac{s}{2})\,.$$

We have

$$\phi(t) = 2e^{-t^2}, \qquad \Theta(t) = \sum_{n=1}^{\infty} 2e^{-\pi n^2 t^2}\,.$$

The function $L^*(s)$ has simple poles in $p_1=0$ and $p_2=1$ with residues $r_1=1, r_2=-1$, so the functional equation for $\Theta(t)$ reads

$$\Theta(1/t) \;=\; t\,\Theta(t) - 1 + t\,. \qquad (7)$$

In fact, applying Poisson's summation formula to $f(x)=e^{-\pi x^2}$ gives (7) and this *proves* the functional equation for $\zeta(s)$.

We now proceed to the second problem, that of computing $L(s)$ or $L^{(m)}(s)$. Fix $s \in \mathbf{C}$ and let

$$G_s(t) = t^{-s} \int_t^{\infty} \phi(x)\, x^s\, \frac{dx}{x}\,, \qquad t>0\,. \qquad (8)$$

Thus $t^s G_s(t)$ is the incomplete Mellin transform of $\phi(t)$ and $\lim_{t \to 0} t^s G_s(t) = \gamma(s)$ is the original $\Gamma$-factor. Again the function $G_s(t)$ decays exponentially with $t$ and can be effectively computed numerically (§4,5).

Consider (5) which expresses $L^*(s)$ as the Mellin transform of $\Theta(t)$. Split the integral into two and apply functional equation (6) to the second integral:

$$\begin{aligned}
L^*(s) \;=\; & \int_0^{\infty} \Theta(t) t^s \frac{dt}{t} = \int_1^{\infty} + \int_0^1 = \int_1^{\infty} \Theta(t) t^s \frac{dt}{t} + \int_1^{\infty} \Theta(1/t) t^{-s} \frac{dt}{t} = \\
& \int_1^{\infty} \Theta(t) t^s \frac{dt}{t} + \int_1^{\infty} \epsilon t^w \Theta(t) t^{-s} \frac{dt}{t} - \int_1^{\infty} \sum_j r_j t^{p_j} t^{-s} \frac{dt}{t} = \\
& \int_1^{\infty} \Theta(t) t^s \frac{dt}{t} + \epsilon \int_1^{\infty} \Theta(t) t^{w-s} \frac{dt}{t} + \sum_j \frac{r_j}{p_j - s}\,.
\end{aligned} \qquad (9)$$

By definition of $\Theta(t)$ and that of $G_s(x)$, the first integral can be rewritten,

$$\int_1^\infty \Theta(t) t^s \frac{dt}{t} = \int_1^\infty \sum_{n=1}^\infty a_n \phi(\frac{nt}{A}) t^s \frac{dt}{t} = \sum_{n=1}^\infty a_n \int_1^\infty \phi(\frac{nt}{A}) t^s \frac{dt}{t} =$$
$$\sum_{n=1}^\infty a_n \int_{n/A}^\infty \phi(t) (\frac{At}{n})^s \frac{dt}{t} = \sum_{n=1}^\infty a_n G_s(\frac{n}{A}).$$

The same applies to the second integral (replace $s$ by $w - s$), so (9) becomes

$$L^*(s) = \sum_{n=1}^\infty a_n G_s(\frac{n}{A}) + \epsilon \sum_{n=1}^\infty a_n G_{w-s}(\frac{n}{A}) + \sum_j \frac{r_j}{p_j - s}.$$

This formula allows to determine $L^*(s)$ and hence $L(s) = L^*(s)/\gamma(s)$ for a given $s \in \mathbf{C}$. Differentiation also gives the formula for the derivatives,

$$\frac{\partial^k}{\partial s^k} L^*(s) = \sum_{n=1}^\infty a_n \frac{\partial^k}{\partial s^k} G_s(\frac{n}{A}) + \epsilon \sum_{n=1}^\infty a_n \frac{\partial^k}{\partial s^k} G_{w-s}(\frac{n}{A}) + \sum_j \frac{r_j (k-1)!}{(p_j - s)^k} \tag{10}$$

It remains to explain how to compute the functions $\phi(t)$ and $\frac{\partial^k}{\partial s^k} G_s(t)$. This occupies the next three sections.

**Remark 2.5.** We assumed that the functional equation (2) involves $L^*(s)$ both on the left-hand and on the right-hand side. In fact, for general motives the functional equation may be of the form

$$L^*(s) = \epsilon \widehat{L}^*(w-s)$$

where

$$L(s) = \sum_{n=1}^\infty \frac{a_n}{n^s}, \qquad \widehat{L}(s) = \sum_{n=1}^\infty \frac{\widehat{a}_n}{n^s}$$

are $L$-functions of "dual" motives. For instance, Dirichlet $L$-series associated to non-quadratic characters are of this nature. Clearly our arguments go through in this more general case as well. The result is that (6) and (10) have to be simply replaced by

$$\Theta(1/t) = \epsilon t^w \widehat{\Theta}(t) - \sum_j \widehat{r}_j t^{p_j}.$$

and

$$\frac{\partial^k}{\partial s^k} L^*(s) = \sum_{n=1}^\infty a_n \frac{\partial^k}{\partial s^k} G_s(\frac{n}{A}) + \epsilon \sum_{n=1}^\infty \widehat{a}_n \frac{\partial^k}{\partial s^k} \widehat{G}_{w-s}(\frac{n}{\widehat{A}}) + \sum_j \frac{\widehat{r}_j (k-1)!}{(\widehat{p}_j - s)^k}$$

Here $\widehat{A}, \widehat{p}_j$ etc. are associated to $\widehat{L}(s)$ as $A, p_j$ etc. are to $L(s)$.

# 3 Computing $\phi(t)$ for $t$ small

Recall that

$$\gamma(s) = \Gamma\Big(\frac{s+\lambda_1}{2}\Big)\cdots\Gamma\Big(\frac{s+\lambda_d}{2}\Big) \tag{11}$$

and $\phi(t)$ is defined as the inverse Mellin transform of $\gamma(s)$. By Mellin's inversion formula (see e.g. [4], §2), $\phi(t)$ is given by a residue sum,

$$\phi(t) = \sum_{z\in\mathbf{C}} \mathrm{res}_{s=z}\, \gamma(s)t^{-s}, \qquad t>0. \tag{12}$$

Since $\Gamma(s)$ has simple poles at $0$ and the negative integers, the function $\gamma(s)$ has a pole at $s\in\mathbf{C}$ iff $s = -\lambda_j - 2n$ for some $j$ and an integer $n$. If $\lambda_j - \lambda_k \notin 2\mathbf{Z}$ for $j\neq k$, then all poles of $\gamma(s)$ are simple and

$$\mathrm{res}_{s=-\lambda_j-2n}(\gamma(s)t^{-s}) = 2\frac{(-1)^n}{n!}\, t^{\lambda_j+2n} \prod_{k\neq j} \gamma\big(\tfrac{(-\lambda_j-2n)+\lambda_k}{2}\big).$$

Hence in this case (12) is of the form $\sum_j t^{\lambda_j} p_j(t^2)$ where $p_j(t)$ a power series in $t$. The coefficients of $p_j(t)$ satisfy a simple linear recursion coming from the relation $\Gamma(s+1)=s\Gamma(s)$.

**Example 3.1.** Let $d=1$ and let $\lambda_1$ be arbitrary. Then $\phi(t)$ is given by

$$\phi(t) = t^{\lambda_1}\sum_{n=0}^{\infty} 2\frac{(-1)^n}{n!} t^{2n} = 2t^{\lambda_1} e^{-t^2}.$$

In general, the poles of $\gamma(s)$ are not simple and the residue of $\gamma(s)t^{-s}$ in $s=z$ is $t^{-z}$ times a polynomial in $\ln t$ of a fixed degree. The reason is that non-constant terms of the Taylor expansion of $t^{-s}$ around $s=z$,

$$t^{-s} = t^{-z}\sum_{k=0}^{\infty} \frac{(-\ln t)^k}{k!}(s-z)^k$$

contribute to the residue in case of a multiple pole. So (12) is again of the form $\sum_j t^{\lambda_j}p_j(t^2)$, except now $p_j(t)$ is a power series in $t$ whose coefficients are polynomials in $\ln t$ of a fixed degree depending on $j$.

**Example 3.2.** Let $d=2$ and $\lambda_1 = \lambda_2 = 0$. Then $\phi(t)$ is a Bessel function,

$$\phi(t) = 4K_0(2t) = -4(\ln t + \gamma_e) - 4(\ln t - 1 + \gamma_e)t^2 - \tfrac{2\ln t - 3 + 2\gamma_e}{2}t^4 + \ldots$$

where $\gamma_e = -\Gamma'(1)$ is the Euler constant.

**Algorithm 3.3. Expansion of $\phi(t)$ for $t$ small.** The following describes the recursions necessary to determine the coefficients of (12) for a general $\Gamma$-factor $\gamma(s)$.

1. Let $\gamma(s)$ and $\phi(t)$ be defined by (11) and (12) respectively.

2. We say that $\lambda_j$ and $\lambda_k$ are equivalent if $\lambda_j - \lambda_k \in 2\mathbf{Z}$. Let $\Lambda_1, ..., \Lambda_N$ denote the equivalence classes and let $l_j = |\Lambda_j|$. Thus $\sum l_j = d$.

3. Let $m_j = -\lambda_{k_j} + 2$ where $\lambda_{k_j} \in \Lambda_j$ is the element with the smallest real part, that is $\inf_{\lambda \in \Lambda_j} \operatorname{Re} \lambda = \operatorname{Re} \lambda_{k_j}$. In other words, $\gamma(s)$ is analytic at $s = m_j$, has a pole (of some order) at $s = m_j - 2$ and a pole of order $l_j$ at $s = m_j - 2n$ for $n \gg 0$.

4. Let $c_j^{(0)}(s)$ be the beginning of the Taylor series of $\gamma(s + m_j)$ around $s = 0$ which ends with $O(s^{l_j})$ as the last term.

5. For $1 \le j \le d$ and $n \ge 1$ define $c_j^{(n)}(s)$ recursively by

$$c_j^{(n)}(s) = c_j^{(n-1)}(s) / \prod_{k=1}^{d} \left( \tfrac{s+\lambda_k+m_j}{2} - n \right), \tag{13}$$

considered as a quotient of Laurent series in $s = 0$. Note that $c_j^{(n)}(s)$ ends with $O(1)$ for $n \gg 0$. Let $c_{j,k}^{(n)}$ denote the coefficient of $s^{-k}$ in $c_j^{(n)}(s)$.

6. For $t$ real positive, $\phi(t)$ is given by

$$\phi(t) = \sum_{j=1}^{N} t^{-m_j} \sum_{n=1}^{\infty} \Bigl( \sum_{k=0}^{l_j-1} \tfrac{(-\ln t)^k}{k!} c_{j,k+1}^{(n)} \Bigr) t^{2n}. \tag{14}$$

**Remark 3.4.** The series above converges exponentially fast since

$$\max_{j \le N, k \le l_j} |c_{j,-k}^{(n)}| = O((n!)^{-d}), \qquad \text{as } n \to \infty.$$

However, for large $t$ this is not a very effective way to compute $\phi(t)$. Take for instance the series $e^{-t} = \sum_{n=0}^{\infty} (-t)^n / n!$ for $t = 50$. The terms grow up to $3. \times 10^{20}$ for $n = 50$ before starting to decrease to 0 in absolute value. Thus to determine $e^{-50}$ to 10 decimal digits with this series one has to require working precision of 30 digits and compute 160 terms until everything happily cancels leading to the answer 0.0000000000. This is clearly not terribly effective. As this is exactly the general behaviour for arbitrary $\gamma(s)$ we use a different method for large $t$, based on the asymptotic expansions at infinity. This is described in §5 below.

## 4 Computing $G_s(t)$ for $t$ small

As explained in §2, we also need a way to compute the incomplete Mellin transform of $\phi(t)$ and its derivatives. Recall that for $s \in \mathbf{C}$ and $t > 0$ we define $G_s(t)$ to be

$$G_s(t) = t^{-s} \int_t^{\infty} \phi(x) x^s \frac{dx}{x}.$$

Recall also that $\lim_{t \to 0} t^{-s} G_s(t)$ exists and equals $\gamma(s)$ whenever $s$ is not a pole of $\gamma(s)$. For such $s$ clearly

$$t^s G_s(t) = \gamma(s) - \int_0^t \phi(x) x^s \frac{dx}{x}. \tag{15}$$

Since (14) expresses $\phi(t)$ as an infinite sum of terms of the form $t^\alpha (\ln t)^\beta$, term-wise integration in (15) gives a similar expression for $G_s(t)$.

In the points where $\gamma(s)$ does have a pole, the formula (15) does not make sense as the right-hand side becomes $\infty - \infty$. However, it is not difficult to locate the terms which give contributions to the principal parts of the Laurent series. Then ignoring these terms gives the value of $G_s(t)$ for such $s$.

**Algorithm 4.1. Expansion of $\frac{\partial^k}{\partial s^k} G_s(t)$ for $t$ small.** This is all summarised in the following formulae which allow to determine $\frac{\partial^k}{\partial s^k} G_s(t)$ for arbitrary $s \in \mathbf{C}$ and $t > 0$. Here $\alpha \in \mathbf{C}$ and $i, j, k \geq 0, n \geq 1$ are integers.

1. Let $c_{j,i}^{(n)}$ be as in (13).
2. Define $L_{\alpha,j,k}(x) \in \mathbf{C}[x]$ by the formula

$$L_{\alpha,j,k}(x) = \begin{cases} k! \sum_{i=0}^{j-1} \binom{i-j}{k} \frac{\alpha^{i-j-k}}{i!} (-x)^i, & \alpha \neq 0, \\ 0, & \alpha = 0. \end{cases}$$

3. Let

$$S_{j,k,s}^{(n)}(x) = \sum_{i=1}^{l_j} c_{j,i}^{(n)} L_{2n+s-m_j,i,k}(x) \quad \in \mathbf{C}[x].$$

4. For $t > 0$ consider the infinite sum

$$\tilde{G}_{s,k}(t) = \sum_{j=1}^{N} t^{2-m_j} \sum_{n=1}^{\infty} S_{j,k,s}^{(n)}(\ln t) \; t^{2n} \tag{16}$$

5. The formula for $\frac{\partial^k}{\partial s^k} G_s(t)$ is

$$\frac{\partial^k}{\partial s^k} G_s(t) = \left( \frac{\partial^k}{\partial S^k} \frac{\gamma(S)}{t^S} \right)^*_{S=s} - \tilde{G}_{s,k}(t),$$

where $f(S)^*_{S=s}$ denotes the constant term $a_0$ of the Laurent expansion $\sum_k a_k (S-s)^k$ of $f(S)$ at $S=s$. Thus $f(S)^*_{S=s} = f(s)$ if $f(S)$ is analytic at $S=s$.

**Remark 4.2.** The series for $\frac{\partial^k}{\partial s^k} G_s(t)$ converges exponentially fast since the corresponding one for $\phi(t)$ does (cf. 3.4). Again, however, it is not effective for large $t$, in which case we use an alternative approach described in the following section.

## 5 Computing $\phi(t)$ and $G_s(t)$ for $t$ large

To compute $\phi(t)$ and $G_s(t)$ for large $t$ we begin with the asymptotic expansions for these functions around infinity.

Recall that $\phi(t)$ is defined as the inverse Mellin transform of a product of $\Gamma$-functions,

$$\Gamma\left(\frac{s+\lambda_1}{2}\right) \cdots \Gamma\left(\frac{s+\lambda_d}{2}\right) = \int_0^\infty \phi(t) \, t^s \, \frac{dt}{t}.$$

In other words, $\phi(t)$ is a special case of Meijer $G$-function. Given two sequences of complex parameters,

$$a_1, \ldots, a_n, a_{n+1}, \ldots a_p \quad \text{and} \quad b_1, \ldots, b_m, b_{m+1}, \ldots b_q$$

a general Meijer $G$-function $G_{p,q}^{m,n}(t; a_1, ..., a_p; b_1, ..., b_q)$ is defined by

$$\int_0^\infty G_{p,q}^{m,n}(t; a_1, ..., a_p; b_1, ..., b_q) t^s \frac{dt}{t} = \frac{\prod_{j=1}^m \Gamma(b_j+s) \prod_{j=1}^n \Gamma(1-a_j-s)}{\prod_{j=m+1}^q \Gamma(1-b_j-s) \prod_{j=n+1}^p \Gamma(a_j+s)}$$

We refer to Luke [17], 5.2-5.11 for basic properties of the $G$-function.

In our case replacing $s$ by $s/2$ yields an identification

$$\phi(t) = 2\, G_{0,d}^{d,0}(t^2; ; \frac{\lambda_j}{2})$$

As discovered by Meijer (in greater generality), the function $G_{0,d}^{d,0}$ possesses the following asymptotic expansion at infinity ([17], Theorem 5.7.5)

$$G_{0,d}^{d,0}(t; ; \frac{\lambda_j}{2}) \sim \frac{(2\pi)^{(d-1)/2}}{\sqrt{d}} e^{-d\, t^{1/d}} t^{\kappa/d} \sum_{n=0}^\infty M_n\, t^{-n/d}$$

$$\kappa = (1 - d + \sum_{j=1}^d \lambda_j)/2 \tag{17}$$

Here $M_n = M_n(\lambda_1, ..., \lambda_d)$ are constants, $M_0 = 1$. As for $\phi(t)$, it follows that

$$\phi(t^{d/2}) \sim \frac{2(2\pi)^{(d-1)/2}}{\sqrt{d}} e^{-d\, t} t^\kappa \sum_{n=0}^\infty M_n\, t^{-n} \tag{18}$$

We would like to note here that the stated asymptotic expansion for large $t$ is much "neater" than the expansion (14) of $\phi(t)$ for small $t$: it does not involve any logarithmic terms and its shape is independent of whether any of the $\lambda_j$ are equal modulo $2\mathbf{Z}$.

The coefficients $M_n$ in the asymptotic expansion can be found as follows. The defining relation

$$\gamma(s+2) = \gamma(s) \times \prod_{j=1}^d \frac{s+\lambda_j}{2}$$

on the level of inverse Mellin transforms is equivalent to an ordinary differential equation (of degree $d$) with polynomial coefficients for $\phi(t)$. It follows that the function $t^{-\kappa} e^{dt} \phi(t^{d/2})$ satisfies a different ODE, of degree $d+1$. Formally substituting $1 + \sum_{n\geq 1} M_n t^{-n}$ as a solution gives a recursion for the $M_n$ with polynomial coefficients. This has been worked out in general by E. M. Wright; see Luke [17], 5.11.5, especially formulae (8) and (16) for details.

**Algorithm 5.1. Asymptotic expansion associated to $\phi(t)$.** Here is the answer in our case, re-written in a slightly different polynomial basis.

1. Let $S_m = S_m(\lambda_1, ..., \lambda_d)$ denote the $m$-th elementary function of $\lambda_1, ..., \lambda_d$,

$$S_0 = 1, \quad S_1 = \sum_{j=1}^{d} \lambda_j, \quad \ldots, \quad S_d = \prod_{j=1}^{d} \lambda_j.$$

2. Define also modified symmetric functions $\tilde{S}_{d+1} \equiv 0$ and

$$\tilde{S}_m = \sum_{k=0}^{m} (-S_1)^k d^{m-1-k} \binom{k+d-m}{k} S_{m-k}, \quad 0 \leq m \leq d.$$

3. For $k \geq 0$ construct $\Delta_k(x) \in \mathbf{Q}[x]$ by means of the generating function

$$\left(\frac{\sinh t}{t}\right)^x = \sum_{k=0}^{\infty} \Delta_k(x) t^{2k}$$

4. For $p \geq 1$ consider the following polynomials

$$\nu_p(n) = -\frac{d}{(2d)^p} \sum_{m=0}^{p} \tilde{S}_m \prod_{j=m}^{p-1} (d-j) \sum_{k=0}^{\lfloor \frac{p-m}{2} \rfloor} \frac{(2n-p+1)^{p-m-2k}}{(p-m-2k)!} \Delta_k(d-p)$$

5. The coefficients $M_n$ in the asymptotic expansion (18) satisfy a recursion

$$M_n = \begin{cases} 0, & n < 0, \\ 1, & n = 0, \\ \frac{1}{n} \sum_{p=1}^{d} \nu_p(n) M_{n-p}, & n \geq 1. \end{cases}$$

Applying term-wise integration to (18), it is also easy to deduce the asymptotic expansion for $G_s(t)$ for $t \to \infty$,

$$G_s(t^{d/2}) \sim \frac{(2\pi)^{(d-1)/2}}{\sqrt{d}} e^{-dt} t^{\kappa-1} \sum_{n=0}^{\infty} \mu_n(s) t^{-n}. \tag{19}$$

Here $\kappa = (1 - d + S_1)/2$ as in (17) and $\mu_n(s) = \mu_n(\lambda_1, ..., \lambda_d; s)$ satisfy a recursion

$$\mu_n = \begin{cases} 0, & n < 0, \\ 1, & n = 0, \\ \frac{1}{n} \sum_{p=1}^{d} \left(\nu_{p+1}(n) - \frac{S_1 + d(s-1) - 2(n-p) - 1}{2d} \nu_p(n)\right) \mu_{n-p}, & n \geq 1. \end{cases} \tag{20}$$

By induction one shows that $\mu_n$ is a polynomial in $s$ with the leading term $2^{-n} s^n$. So if we differentiate (19) $k$ times to $s$, exactly $k$ terms vanish and we get the following formula for the derivatives $\frac{\partial^k}{\partial s^k} G_s(t)$,

$$\frac{\partial^k}{\partial s^k} G_s(t^{d/2}) \sim \frac{(2\pi)^{(d-1)/2}}{\sqrt{d}} e^{-dt} t^{\kappa-1-k} \sum_{n=0}^{\infty} \frac{\partial^k \mu_{n+k}(s)}{\partial s^k} t^{-n}. \tag{21}$$

Equations (18), (19) and (21) provide asymptotic series for the functions $\phi(t)$, $G_s(t)$ and $\frac{\partial^k}{\partial s^k} G_s(t)$ at infinity. For computational purposes, though, it is better to work with continued fraction expansions associated to these series. Consider, for instance, the case of $\phi(t)$, the case of $\frac{\partial^k}{\partial s^k} G_s(t)$ being analogous.

Fix $d$ and $\lambda_1, ..., \lambda_d$. Letting $x = 1/t$ in (18) we get

$$\psi(x) := \frac{\sqrt{d}}{2(2\pi)^{(d-1)/2}} e^{-dx} x^{\kappa} \phi(x^{-d/2}) \sim \sum_{n=0}^{\infty} M_n x^n \tag{22}$$

with $M_n$ constants. As any formal series, the right-hand side can be *formally* written either as a unique infinite continued fraction

$$\sum_{n=0}^{\infty} M_n x^n = \alpha_0 + \cfrac{x^{k_0}}{\alpha_1 + \cfrac{x^{k_1}}{\alpha_2 + \cfrac{x^{k_2}}{\alpha_3 + ...}}}, \qquad \alpha_n \neq 0 \text{ for } n > 0,$$

or as a unique terminating fraction of the same form. To see this start with $p_0(x) = \sum M_n x^n$ and define recursively formal power series $p_{n+1}(x)$ in terms of $p_n(x)$ by

$$p_n(x) = \alpha_n + \frac{x^{k_n}}{p_{n+1}(x)}, \qquad n \geq 0.$$

Here $k_n$ is the degree of the first non-zero term in $p_n(x) - p_n(0)$; if $p_n(x) \equiv 0$ for some $n$, then terminate. This shows the existence of the continued fraction expansion and uniqueness is not difficult to verify as well. This construction also gives an explicit way to calculate the $\alpha_n$ if the $M_n$ are given (up to some limit). There are of course better (computationally more stable) methods, see for instance [11, 16].

If the fraction does not terminate, define the partial convergents $C_n(x)$ for all $n$ by

$$C_n(x) = \alpha_0 + \cfrac{x^{k_0}}{\alpha_1 + \cfrac{x^{k_1}}{... + \cfrac{x^{k_{n-1}}}{\alpha_n}}}.$$

If the fraction does terminate at $C_N$, let $C_n = C_N$ for $n > N$.

In any case we can think of $C_n(x)$ as approximants to the original function $\psi(x)$ of (22). Indeed, $C_n(x)$ is a rational function whose Taylor expansion around $x=0$ starts at least with $M_0 + ... + M_n x^n$. Hence $\psi(x)$ and $C_n(x)$ have the same asymptotic expansions up to $x^n$. Thus there is a constant $K_n > 0$ such that

$$|\psi(x) - C_n(x)| \leq K_n x^{n+1}.$$

Unfortunately, it seems very difficult to provide explicit bounds for $K_n$. It appears that $C_n(x)$ converge rapidly to $\psi(x)$ but to prove either "converge" or "rapidly" or "to $\psi(x)$" in any generality seems hard. So the last step of the algorithm is based purely on empirical observations concerning the convergence of the continued fractions. If one is uncomfortable with this, see §6 on how to avoid this. In the implementation [7] we do use asymptotic expansions with a simple numerical check (see step 7 below) to justify the values.

**Algorithm 5.2. Computing $\phi(t)$ for $t$ arbitrary.** The computation of $\phi(t)$ for arbitrary $t$ can be performed as follows:

1. Let $\epsilon > 0$ be the necessary upper bound for the required precision in the computations of $\phi(t)$.
2. Let $\phi_n(t)$ be the n-th approximant to $\phi(t)$ defined by (cf. (18))

$$\phi_n(t) = \frac{2(2\pi)^{(d-1)/2}}{\sqrt{d}} e^{-d\,t^{2/d}} t^{2\kappa/d} C_n(1/t^{2/d}), \qquad n \geq 0.$$

As we already noted, $\phi(t) - \phi_n(t) = O(t^{-n})$ as $t \to \infty$. Denote by $\phi_{taylor}(t)$ the function $\phi(t)$ computed using the power series expansion at the origin as in §3.

3. Determine $t_0$ such that $|\phi_0(t)| < \epsilon/2$ for $t > t_0$.
4. Choose a subdivision of the interval $[0, t_0]$,

$$0 < t_k < t_{k-1} < \ldots < t_1 < t_0 < \infty$$

For every $t_i$ let $n_i$ be an integer for which $|\phi_n(t) - \phi_{n+1}(t)| < \epsilon/2$ and $|\phi_n(t) - \phi_{n+2}(t)| < \epsilon/2$.

5. Determine $M_n$ for $0 \leq n \leq n_k$ using the recursion they satisfy.
6. The function $\phi(t)$ is computed as follows

$$\phi_{general}(t) = \begin{cases} \phi_{taylor}(t), & 0 < t \leq t_k \\ \phi_{n_i}(t), & t_i \leq t < t_{i-1} \\ 0, & t > t_0 \end{cases}$$

7. As a numerical check, verify that $\phi_{taylor}(t_k)$ agrees with $\phi_{n_k}(t_k)$.

**Example 5.3.** Let $d = 2$ and $\lambda_1 = \lambda_2 = 0$ as in Example 3.2. Recall that $\phi(t) = 4K_0(2t)$ is a Bessel function in this case. Asymptotic expansion (18) then reads

$$\phi(t) \sim 2\sqrt{\pi}\, e^{-2t}\, t^{-1/2} \sum_{n=0}^{\infty} M_n\, t^{-n}$$

and the coefficients $M_n$ satisfy a recursion

$$16 n M_n = -(2n-1)^2 M_{n-1}.$$

It follows that

$$M_0 = 1, \ M_1 = -\tfrac{1}{16}, \ M_2 = -\tfrac{9}{512}, \ \ldots, M_n = -\frac{(2n-1)!!(2n-1)!!}{16^n n!}, \ \ldots$$

Choose $\epsilon = \tfrac{1}{2} \cdot 10^{-10}$ and $t_0 = 12$, $t_1 = 6$, $t_2 = 2$. Take $n_1 = 6$ and $n_2 = 20$ and compute $\phi(t)$ by

$$\phi_{general}(t) = \begin{cases} \phi_{taylor}(t), & 0 < t \leq 2 \\ \phi_{20}(t), & 2 \leq t < 7 \\ \phi_6(t), & 7 \leq t < 12 \\ 0, & t > 12 \end{cases}$$

As a numerical check we verify that $|\phi_{taylor}(2) - \phi_{20}(2)| \leq 4 \cdot 10^{-14} \leq \epsilon$ as required.

## 6 Implementation notes

Let us begin with a summary of steps of the algorithm presented in the previous sections. We start with an $L$-function satisfying 2.1–2.3 (see also Remark 2.5).

- The formula used for the numerical verification of the functional equation is (6) and that for computing $L(s)$ and its derivatives is (10) together with $L^*(s) = L(s)/\gamma(s)$. The functions used in these formulae are $\phi(t)$ defined by (3) and $G_s(t)$ defined by (8).

- To compute $\phi(t)$ numerically we employ (5.2). It is based on power series expansions in the origin (3.3), asymptotic formula (18), recursion (5.1) and the associated continued fractions.

- Similarly $\frac{\partial^k}{\partial s^k} G_s(t)$ is computed in the same way. The corresponding expansion in the origin is given by (4.1), asymptotics by (19) and recursion for the coefficients by (20).

However, in order to make a practical algorithm out of these results, we still need to explain how to truncate various infinite sums and discuss related precision issues.

If one desires to implement our method with rigorous proofs that all of the computations are correct, the following issues have to be considered. First, one has to keep track of the number of operations used and the possible round-off errors, perhaps even using interval arithmetic to justify the computations. Second, one needs to have analytic bounds on the size of the functions $\phi(t)$ and $G_s(t)$ for large $t$, rather than just asymptotic behaviour.

In the PARI implementation [7] we have chosen to be content with the intuitively natural bounds and a few numerical checks to justify the results. A reader wishing to use a more rigorous approach, might consider the following:

**Remark 6.1.** Let us start with the computations of $\phi(t)$ and $\frac{\partial^k}{\partial s^k} G_s(t)$ by means of series expansions around the origin. Both (14) and (16) defining these are infinite sums, but it is not difficult to see how to terminate them. The point is that it suffices to give an explicit bound on the coefficients $c_{i,j}^{(n)}$ which goes to 0 exponentially with $n$. Everything else in (14) and (16) grows at most polynomially in $n$ so any rough estimate on them will do. As for an explicit exponential bound on $c_{i,j}^{(n)}$, it follows from (13) and, say, the obvious lower bound $\frac{1}{2} n^d$ for $n \gg 0$ on the coefficients in $n^d \prod_{k=1}^d (1 - \frac{s+\lambda_k+m_j}{2n})$ treated as a polynomial in $s$.

**Remark 6.2.** The next question is that of working precision. This has already been mentioned in 3.4 and 4.2. When $\phi(t)$ and $G_s(t)$ are computed from the expansions around the origin, the terms in the defining series can be very large. Therefore one needs to work with larger working precision than the desired precision of the answer.

Similar cancellation in fact occurs when one computes $L(s)$ for $s$ with large imaginary part. This is a well-known problem, which one has to face when verifying the Riemann hypothesis. The point is that $L(s) = L^*(s)/\gamma(s)$ and both $L^*(s)$ and $\gamma(s)$ decrease exponentially fast on the vertical strips. Hence one needs to compute $L^*(s)$ to more significant digits ($\log_{10} |\gamma(s)|$ more to be precise), to evaluate $L(s)$ to given precision.

A solution to this has been suggested in Lagarias-Odlyzko [15] and worked out by Rubinstein [19], at least for $d=1$ and $d=2$. By modifying $G_s(t)$ by a suitably chosen exponential factor one obtains a formula for $L(s)$ which does not have the loss-of-precision behaviour. It is perhaps possible to use Rubinstein's approach to obtain similar formulae for a general $\Gamma$−factor as well.

**Remark 6.3.** The next issue is that of truncating the main formulae used is this paper. Recall that to verify the functional equation numerically, we use the function

$$\Theta(t) = \sum_{n=1}^{\infty} a_n \, \phi(\tfrac{nt}{A}) \,. \tag{23}$$

Then to actually compute the $L$-values, we use

$$\frac{\partial^k}{\partial s^k} L^*(s) = \sum_{n=1}^{\infty} a_n \frac{\partial^k}{\partial s^k} G_s(\tfrac{n}{A}) + \epsilon \sum_{n=1}^{\infty} a_n \frac{\partial^k}{\partial s^k} G_{w-s}(\tfrac{n}{A}) + \sum_j \tfrac{r_j(k-1)!}{(p_j-s)^k} \,. \tag{24}$$

See also Remark 2.5 for the necessary modifications when there are two different $L$-functions involved in the functional equation. In any case, one needs analytic estimates on $\phi(t)$ and $\frac{\partial^k}{\partial s^k} G_s(t)$ for large $t$ to carefully estimate the error in truncating these infinite sums.

One possible way to obtain such estimates is to use Tollis' method [23] based on Braaksma's work [4] on asymptotic behaviour of Meijer G-functions. By applying the Euler-Maclaurin summation formula to the Mellin-Barnes integral defining $G_s(t)$, Tollis determines an explicit exponential bound for $G_s(t)$ in the case $\lambda = (0,...,0,1,...,1)$ with $\rho+\sigma$ zeroes and $\sigma$ ones ($\rho, \sigma \geq 0$). It is likely that his method is general enough to obtain similar estimates for an arbitrary $\Gamma$-factor as well.

**Remark 6.4.** Finally, let us turn to the methods of §5, asymptotic expansions and associated continued fractions.

Unfortunately, there seem to be few cases where one can actually provide explicit estimates for the convergence of the continued fractions of, say, $\phi(t)$. The most general result known to the author in this connection is that of Gargantini and Henrici [10]. They show that functions which can be written as Stieltjes transforms of positive measures admit convergent continued fraction expansions at infinity, with explicit error bounds. This does not seem to apply to our functions in general, though. See Henrici [11], Chapter 12 and Lorentzen-Waadeland [16] for more information.

The analysis is available, though, in low-dimensional cases. For instance for $d=1$ the function $G_s(x)$ is the incomplete Gamma function, for which there are known convergent continued fractions expansions at infinity, see Henrici [11], 12.13.I. Also for $d=2$ the function $\phi(x)$ reduces to the Whittaker function, which is a Stieltjes transform (basically, of itself). So in this case the continued fraction expansion converges, see Henrici [11], 12.13.II.

One possible way out is to compute $\phi(t)$ and $\frac{\partial^k}{\partial s^k} G_s(t)$ only using Taylor expansions at the origin, even for large $t$. It is easier to give precise estimates for the convergence in this case, although one does pay the price with substantial loss of efficiency. Alternatively, one might try a completely different approach to compute the functions in question at infinity. For instance, it is perhaps possible to use backward recursions, as one does for Bessel functions.

# 7 *L*-functions with unknown invariants

In a perfect world, one knows all of the invariants associated to one's *L*-function. In a less perfect world, one does not know exactly the sign $\epsilon$ and, perhaps, the residues $r_i$ at the poles of $L^*(s)$. In reality, however, there are plenty of examples where it is difficult to determine the exponential factor $A$ and even some of the coefficients $a_n$. Fortunately, in some of these cases it is still possible to make computations with *L*-functions.

To illustrate this, say that $L(s)$ is expected to satisfy Assumptions 2.1-2.2 and only the sign $\epsilon$ in the functional equation is difficult to determine. As we already mentioned, the functional equation is equivalent to the statement that for all $1<t<\infty$,

$$\Theta(1/t) = \epsilon\, t^w \Theta(t) - \sum_j r_j t^{p_j}. \tag{25}$$

Choose $1<t<\infty$ and evaluate the left-hand and the right-hand side. This gives an equation which can be solved for $\epsilon$. Of course it is then sensible to verify that (25) holds with the obtained $\epsilon$ by verifying it numerically for some other values of $t$.

The same applies in the case when neither $\epsilon$ nor the residues $r_i$ are known. The equation above is linear in all of these, so choosing enough $t$'s gives a linear system of equations from which $\epsilon$ and the $r_i$ can be obtained. There might be of course precision problems if there are many residues to be determined.

In most cases, actually, $\epsilon = \pm 1$ and $L^*(s)$ has no poles, so simply trying $\epsilon = -1$ and $\epsilon = 1$ for some $t>1$ immediately yields the right sign.

Next come the dimension $d$, the Hodge numbers $\lambda_1,...,\lambda_d$ and the poles $p_j$ of $L^*(s)$. Fortunately, these can always be determined in practice, at least in all of the cases that the author is aware of.

The next issue is that of the exponential factor $A$. For instance, consider $L(C, H^1, s)$, the *L*-function associated to $H^1$ of a genus $g$ curve $C/\mathbf{Q}$. Then $A = \sqrt{N}/\pi^g$ where $N$ is the conductor of $C$. In practice, to determine $N$ one needs at least to be able to find a model of $C$ over $\mathbf{Z}$ which is regular at a given prime of bad reduction. This, in turn, means performing successive blowing-ups over the unramified closure of $\mathbf{Q}_p$, an operation not without computational difficulties. For curves of genus 1 and 2 there are effective algorithms for doing this, but not for higher genus. So finding $N$ for a given curve might be hard in practice. Also note that (25) is absolutely *not* linear in $A$, so one cannot solve for it directly.

Fortunately, one can usually determine the full set $\Sigma = \{p_1,...,p_k\}$ of primes where $C$ has bad reduction. Then one knows that $N = p_1^{b_1} \cdots p_k^{b_k}$ is composed of those primes and has (hopefully) an upper bound for the $b_i$, say in terms of the discriminant of $C$ or some similar quantity. This leaves only finitely many choices for $N$ and (as in the case of the sign $\epsilon = \pm 1$), a simple trial-and-error can establish the proper functional equation. It should be noted here that this applies, of course, only to those *L*-functions for which there is a *unique* $A$ (and $\epsilon$ etc.) for which the functional equation holds.

Finally we come to the coefficients $a_i$. Again take the case of a genus $g$ curve $C/\mathbf{Q}$ with the set $\Sigma = \{p_1,...,p_k\}$ of bad primes as above. Then the problem is to determine the local factors at bad primes, that is the coefficients $a_{p_i^j}$ for $1 \leq i \leq k$ and $j \geq 1$. The local factors at good primes can be determined by counting points over finite fields and the coefficients $a_n$ for composite $n$ can be obtained by the product formula.

Fortunately, again, there are only finitely many choices of possible local factors for a given bad prime $p_i$. For instance $|a_p| < 2\sqrt{p}$ and $|a_{p^j}|$ satisfy similar estimates. Moreover, $a_{p^j}$ with $1 \le j < 2g$ determine $a_{p^j}$ for all $j$ as the degree of the local factor is bounded by $2g$. Note, however, that this is not a very practical approach, especially for large primes $p_i$ when there are numerous possibilities for the local factors.

Another approach is to note that the functional equation (25) is in fact *linear* in the $a_i$, since $\Theta(t)$ is. If there were only finitely many unknown coefficients $a_i$, they could obtained in the same way as $\epsilon$ and the $r_i$ were.

To illustrate what can be done when infinitely many coefficients are unknown, consider the following typical case:

1. Say, there is just one prime $p$ for which $a_p$ is difficult to determine theoretically,

2. assume that all $a_n$ are integers,

3. assume that there is a product formula for $L(s)$ in question, in particular $a_{mn} = a_m a_n$ for $m, n$ co-prime.

Using multiplicativity, write (4) as

$$\Theta(t) = \sum_{k=1}^{\infty} a_{p^k} \theta_k(t)$$

where $\theta_k(t)$ are computable functions. Moreover, since we only take finitely many terms when actually computing $\Theta(t)$, we have

$$\Theta(t) \approx \sum_{k=1}^{K} a_{p^k} \theta_k(t)$$

where "$\approx$" stands for "equal to required precision". Hence the functional equation (25) for a fixed $t$ becomes simply a linear equation in $a_p, ..., a_{p^K}$. Thus we can again plug in enough $t$'s to get a linear system which can be solved for the $a_{p^i}$. However, the coefficient functions $\theta_k(t)$ decay rapidly with $k$, so $a_{p^i}$ obtained from solving this system are certainly unreliable for large $i$. If, however, the first coefficient $a_p$ does look like an integer, we can simply round it off and repeat the same process with $a_{p^2}, ..., a_{p^K}$ as variables until all the $a_{p^i}$ are determined.

In practice this works well for a large prime $p$ and even when there are several (large) primes $p$ for which the $a_{p^i}$ are unknown. This does not work for small primes, for instance virtually never for $p = 2$. But then for *small* primes one may try all possible local factors by trial-end-error and for *large* primes solve for the coefficients as described here.

At this moment the reader might be long horrified by the methods suggested here and might wonder whether the reliability of such an approach is not extremely dubious. In our defence we can say that since there is a very effective way to verify the functional equation numerically, *any* method to make an intelligent guess will do, however dubious it might be. When $A$, $\epsilon$ and the bad local factors are determined (or simply guessed in whatever way), one can make numerous checks that these have the correct shape and that (6) holds for various $t$. Thus it is hoped that someone who actually tried to perform blowing-ups on a genus 6 arithmetic surface which has $2^{20}$ in the discriminant will forgive the author for offering desperate tricks to avoid the hard work. After all, this does allow to give evidence for various conjectures even in the difficult cases where it is hard to determine all of the invariants of the $L$-function in question using theoretical arguments.

Finally, let us mention here that there are fortunately better ways to guess the local factors for bad primes, at least for arithmetic surfaces. These have been used to make computations with curves of genus $g \leq 8$ and are to appear in [8].